\documentclass[12pt]{article}

\usepackage{amsmath, amstext, amsgen, amsbsy, amsopn, amsfonts, amssymb, graphicx,psfrag, mdwlist}

\usepackage{pdfsync}
\usepackage{epstopdf}
\usepackage{color}
\usepackage{subfig}

\newtheorem{theorem}{Theorem}

\newtheorem{corollary}[theorem]{Corollary}

\title{A Note on Alexander Polynomials of 2-Bridge Links}
\author{Jim Hoste\\Pitzer College}

\begin{document}
\maketitle
\begin{abstract}
A formula for the Alexander polynomial of a 2-bridge knot or link given by Hartley and also by Minkus has a beautiful interpretation as a walk on the integers. We extend this to the 2-variable Alexander polynomial of a 2-bridge link, obtaining a formula that corresponds to a walk on the 2-dimensional integer lattice.
\end{abstract}
\section{Introduction}
A very elegant  way to compute the Alexander polynomial of a 2-bridge knot (or the reduced Alexander polynomial of a 2-bridge link) was given by Hartley in 1979 \cite{Hartley:1979}. By drawing what Hartley called the ``extended diagram'' of the knot, the Alexander polynomial could simply be read off from the picture. Hartley then used his ``picturesque result'' to prove that the Alexander polynomial of a 2-bridge link is trapezoidal, a property conjectured by Fox (but still unproven in general) to be true for the Alexander polynomial of any alternating knot \cite{Fox:1962}. He also used his approach to provide an algorithm for extremely rapid calculation of the polynomial. A formula similar to Hartley's was given by Minkus in 1982 \cite{Minkus:1982}. Both of these formulae allow the computation of the Alexander polynomial to be beautifully interpreted as a walk on the integers. However, the author has been unable to find such an interpretation anywhere in the literature. In this paper we extend this interpretation to the 2-variable Alexander polynomial of a 2-component, 2-bridge link, where the 1-dimensional walk now becomes a 2-dimensional walk. While it is reasonable to believe that this extension to links may have been known to Hartley or other knot theorists of his era, we are unsure if this is true. Moreover, the  interpretation of the polynomials as walks on integer lattices seem to be little known to current knot theorists. The goal of this short note is to share this viewpoint with a wide audience.

\section{The Alexander Polynomial and 1-Dimensional Walks}

\begin{figure}[!tbp]
  \centering
  \subfloat[Standard diagram of $K_{3/7}$.]{\includegraphics[width=0.4\textwidth]{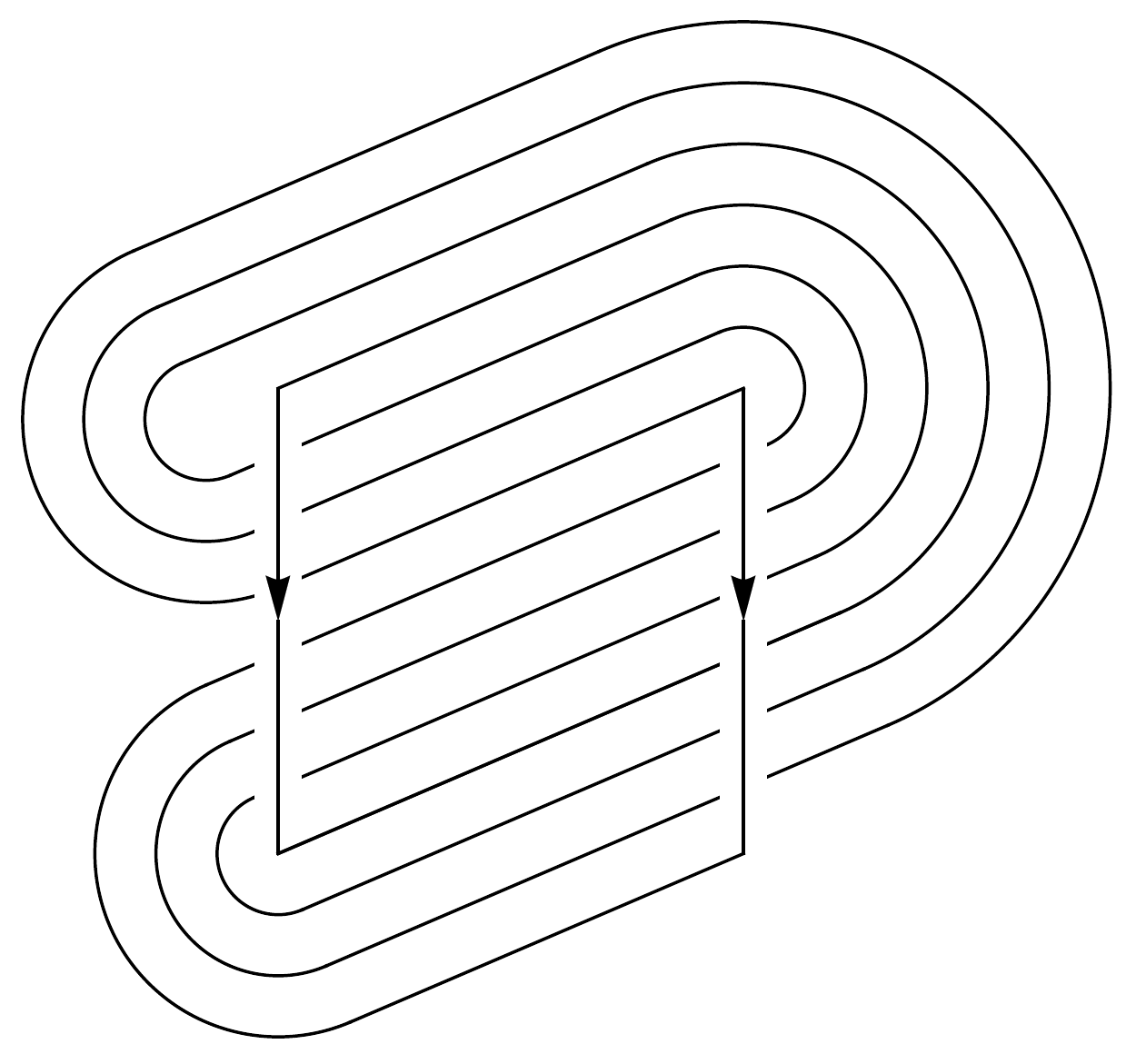}\label{fingerprint}}
  \hfill
  \subfloat[Extended diagram for $K_{3/7}$]{\includegraphics[width=0.4\textwidth]{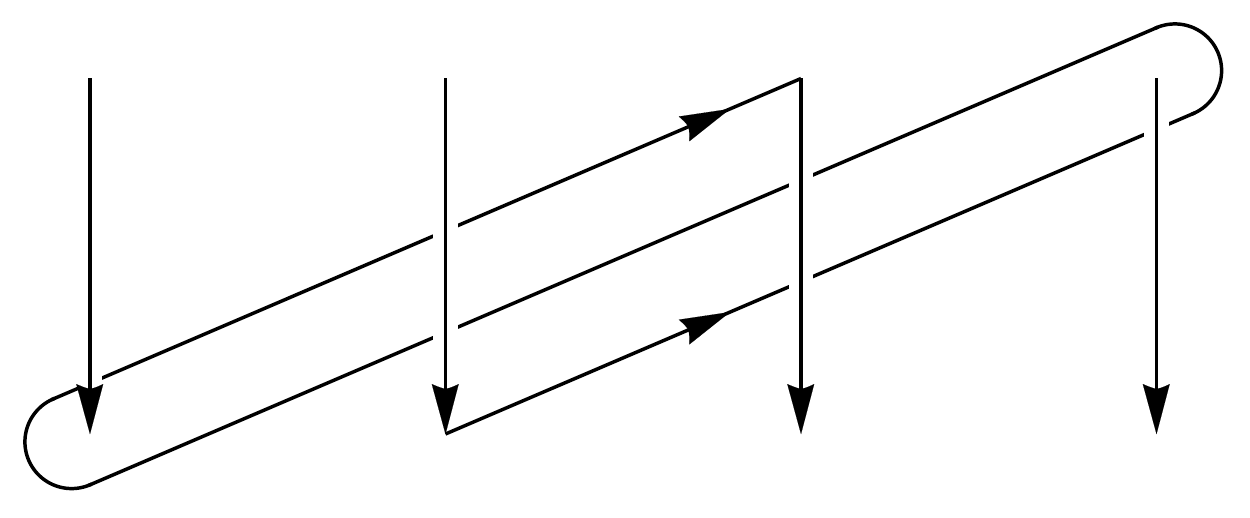}\label{extended}}
  \caption{The standard 2-bridge diagram for $p/q=3/7$ on the left and its extended diagram on the right.}
  \label{std and extended diagrams}
\end{figure}

Recall that each  2-bridge knot or link corresponds to a relatively prime pair of integers $p$ and $q$ with $q$ odd in the case of a knot and $q$ even in the case of a link of two components. We denote the knot or link as $K_{p/q}$. Moreover, we may assume that $0<p<q$ and that $p$ is odd. In this paper, we assume that $K_{p/q}$ is oriented as shown in the {\it standard diagram} of $K_{p/q}$ which appears as in Figure~\ref{fingerprint} (pictured for $p/q=3/7$) where the lines inside the square bounded by the two bridges  have slope $p/q$. If the left bridge is oriented downwards, then, because $p$ is odd, the orientation of the right bridge will also be downwards. If $q$ is even, then each bridge lies in a different component, and the orientation of the bridges as shown in  Figure~\ref{fingerprint} is used to define the orientation of the 2-component link $K_{p/q}$.\footnote{Our use of $p$ and $q$ is reversed from Hartley's.}

The Alexander polynomial $\Delta_K(t)$ of a knot $K$ is only defined up to multiplication by units in $\mathbb Z[t^{\pm1}]$. Therefore, it will be useful to write $f(t)\doteq g(t)$, where $f(t)$ and $g(t)$ are in $\mathbb Z[t^{\pm1}]$, to mean that $f(t) = ug(t)$, where $u$ is a unit in $\mathbb Z[t^{\pm1}]$. The following formula for $\Delta_{K_{p/q}}(t)$ can be found in \cite{Minkus:1982}.

\begin{theorem} [Minkus]  Suppose $K_{p/q}$ is a 2-bridge knot or link with $0 < p < q, p \text{ odd, and}$ $\gcd(p, q) = 1$. Then
\begin{equation} \label{minkus formula}
\Delta_{K_{p/q}}(t)\doteq \sum_{k=0}^{q-1} (-1)^k t^{\sum_{i=0}^k \epsilon_i}
\end{equation}
where $\epsilon_i = (-1)^{\lfloor i p/q\rfloor}$ and $\lfloor x \rfloor$ is the largest integer less than or equal to $x$.
\end{theorem}

We now interpret Minkus' formula for $\Delta_{K_{p/q}}(t)$ as a 1-dimensional walk on the integers, illustrating the algorithm for the case of $p/q = 5/13$. It is not difficult to see that this algorithm produces the polynomial given in  Equation \eqref{minkus formula}.

\begin{figure}[htbp]
\begin{center}
\includegraphics[scale=.6]{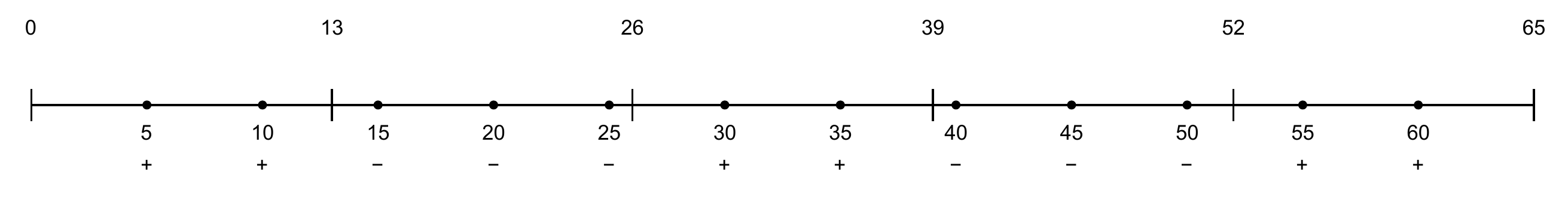}
\caption{The sequence of signs for $p/q=5/13$.}
\label{signs for 5/13}
\end{center}
\end{figure}

\noindent{\bf Algorithm 1:}\footnote{The first author {\sl thought} he learned of this algorithm while visiting the University of Hawaii in the Summer of 1998, presumably from Mike Hilden, but Hilden does not share this recollection!} First mark off the multiples of $q$ from $0$ to $pq$ and the multiples of $p$ from $p$ to $(q-1)p$ as illustrated in Figure~\ref{signs for 5/13} for $p/q = 5/13$. Now write down a sequence of plus and minus signs, one for each multiple of $p$, as follows. For those multiples of $p$ that are less than $q$, record plus signs. Then for multiples of $p$ between $q$ and $2q$, record minus signs. Continue in this fashion, alternately recording plus or minus signs as the multiples of $p$ fall between successive multiples of $q$.

\begin{figure}[htbp]
\begin{center}
\includegraphics[scale=.7]{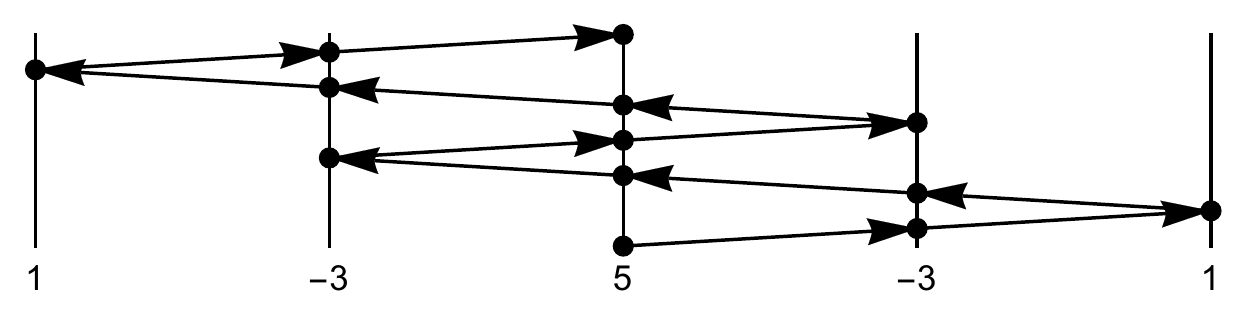}
\caption{The the walk on the integer lattice for $p/q=5/13$ produces the polynomial $1-3t + 5t^2 - 3t^3 + t^4$.}
\label{walk for 5/13}
\end{center}
\end{figure}

We now take a walk on the integers determined by this sequence of signs. Starting at any integer, step to the right for each plus sign and to the left for each minus sign. For $5/13$ the sequence of signs is $++---++---++$. Hence we take two steps to the right, then three to the left, and then two to the right, etc.\,as shown in Figure~\ref{walk for 5/13}. (Because the walk backtracks over itself, we have lifted it up out of the integer lattice for greater clarity.) We now record how many times we visit each integer, obtaining the bi-infinite sequence $\dots, 0,0,1, 3, 5, 3, 1,0,0,\dots$. After ignoring zeroes and introducing an alternation of signs, we obtain the sequence $1, -3, 5, -3, 1$, which are the coefficients of the Alexander polynomial: $\Delta_{K_{5/13}}(t)\doteq 1-3t + 5t^2 - 3t^3 + t^4$.

\noindent{\bf Algorithm 2:} We now describe Hartley's algorithm which is similar to Algorithm 1. We illustrate the algorithm for $p/q = 3/7$. In Figure~\ref{std and extended diagrams} we have pictured the standard 2-bridge diagram in the case of $3/7$ on the left and its associated {\it extended diagram} on the right. To form the extended diagram start with an infinite strip of bridges and then unwind the underarcs as described in \cite{Hartley:1979}. Only one unwound underarc is shown in the figure and for this arc we count the number of segments lying between successive bridges, obtaining the sequence $2, 3, 2$. After introducing alternating signs, these are the coefficients of the Alexander polynomial:  $\Delta_{K_{3/7}}(t) \doteq 2 - 3t + 2t^2$.

Algorithm 2 can also be interpreted as a walk on the integers as determined by the unwound underarc. In the case of $p/q = 3/7$, we take two steps to the right, then three to the left, and then finish with two steps to the right. This sequence of steps can also be derived in almost the same way as Algorithm 1. We continue to illustrate the case of $p/q = 3/7$. First mark off the multiples of $q$ from $0$ to $pq$ as in Figure~\ref{signs for 5/13}. In Algorithm 1 we would then mark off the multiples of $p$. Instead, we now take the multiples of $p$ and shift them down by $p/2$. To this new sequence we again associate plus and minus signs as before, obtaining, in this case, the sequence $+ + --- + +$ which determines the walk. This time, to obtain the coefficients of the polynomial, we count the number of times the walk {\sl passes through each half-integer} and, as before, introduce an alternation of sign. (Alternatively, we could take a walk on the half-integers and then count how many times we pass through each integer.)

\section{Extending to the 2-Variable Polynomial}

In the case of a 2-component, 2-bridge link L, we can compute its 2-variable Alexander polynomial $\Delta_L(x,y)$ in a similar fashion. We illustrate the algorithm in the case of $p/q = 5/18$.

\noindent{\bf Algorithm 3:} As in Algorithm 1, start with the multiples of $p$ and $q$ and produce the corresponding sequence of signs. For the case of $p/q=5/18$, the reader can check that we obtain $+ + + - - - - + + + - - - - + + +$. This sequence of signs determines a walk on the 2-dimensional integer lattice $\mathbb Z \times \mathbb Z$, starting at the point $(0,0)$,  as follows. There is one step in the walk for {\sl every other} sign, starting with the second one. If the sign is positive, the horizontal component of the step is $+1$, otherwise it is $-1$. The vertical component of the step depends on the previous and next signs: if they are different, the vertical component is zero; if they are alike and positive, the vertical component is $+1$ and if they are alike and negative, it is $-1$. Thus in the case of $5/18$, the first step is one unit to the right and one unit up, the second is one unit to the left, etc. If $k$ is the number of times the walk visits the point $(i,j)$, then the Alexander polynomial $\Delta_L(x,y)$ will contain the summand $(-1)^{i+j}k x^i y^j$.

\begin{figure}[htbp]
\begin{center}
\includegraphics[scale=.7]{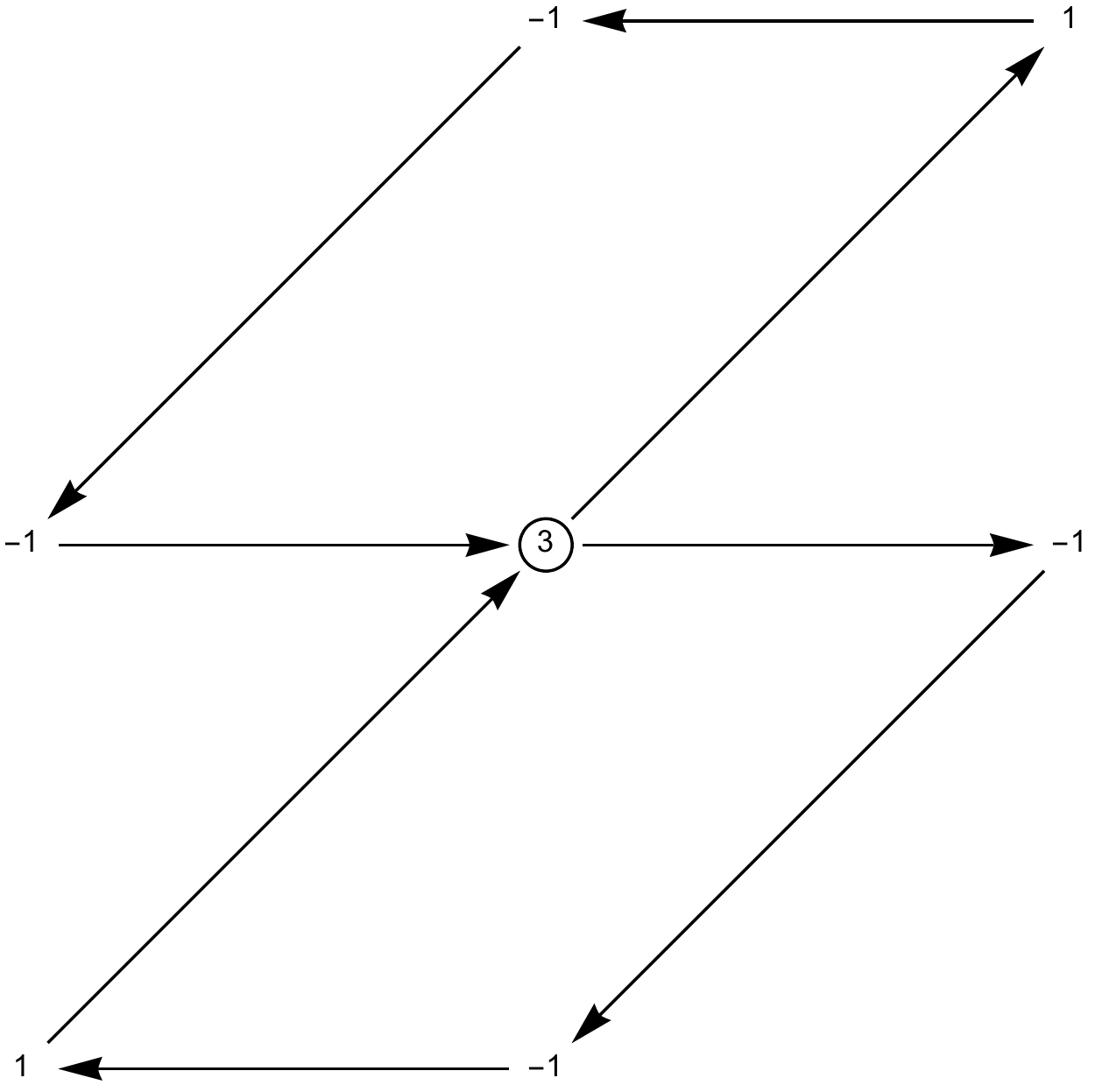}
\caption{The the walk on the 2-dimensional integer lattice for $p/q=5/18$ produces the polynomial $x^2 y^2-x^2 y-x y^2+3 x y-x-y+1$.}
\label{walk for 5/18}
\end{center}
\end{figure}

The walk for $p/q=5/18$ is shown in Figure~\ref{walk for 5/18}. The walk starts at the circled vertex and the coefficient at each vertex is shown. We obtain the polynomial $$\Delta_{K_{5/18}}(x,y)\doteq x^2 y^2-x^2 y-x y^2+3 x y-x-y+1.$$

\section{Derivation of the Algorithms}
As already mentioned, Algorithms 1 and 2 are implicit in the work of Minkus and Hartley. We  sketch a proof of the following theorem and then show that this gives Algorithm 3.

\begin{theorem}\label{2-variable formula}If $K_{p/q}$ is the 2-component,  2-bridge link with $0 < p < q, p \text{ odd}, q \text{ even,}$\\$\text{ and } \gcd(p,q) = 1$, then
\begin{equation}\label{2variableFormula}
\Delta_{K_{p/q}}(x,y)\doteq \sum_{i=1}^{q/2}\epsilon_{2i-1}x^{\sum_{j=1}^{i-1} \epsilon_{2j}}y^{\frac{\epsilon_{2i-1}-1}{2}+\sum_{k=1}^{i-1}\epsilon_{2k-1}}
\end{equation}
where $\epsilon_i=(-1)^{\lfloor \frac{i p}{q} \rfloor}$.
\end{theorem}
\noindent{\bf Proof:}
Starting from the standard 2-bridge diagram of the link $K_{p/q}$ (similar to the one  shown on the left side of Figure~\ref{std and extended diagrams}) we may obtain the following well-known presentation (called the ``over presentation'' in \cite{Hartley:1979}) for the fundamental group $G$ of the complement of $K_{p/q}$:

\begin{equation}\label{group presentation}
G = \langle a,b \, |\, awa^{-1}w^{-1} = 1\rangle
\end{equation}
 where $w = b^{\epsilon_1} a^{\epsilon_2} \dots b^{\epsilon_{q-1}}$ and $\epsilon_i = (-1)^{\lfloor i p/q\rfloor}$. The  signs associated to the
multiples of $p$ in Algorithms 1 and 3 are exactly the signs of the $\epsilon_i$'s.
Following the procedure given in \cite{Crowell_Fox:1963} and using the notation given there, we first use the Fox calculus to compute the Jacobian
 \begin{align*}
 \cal{J}&=\left (\begin{array}{cc} \frac{\partial a w a^{-1}w^{-1}}{\partial a}& \frac{\partial a w a^{-1}w^{-1}}{\partial b} \end{array}\right )\\
 &=\left (\begin{array}{cc} 1-w+(a-1)\frac{\partial w}{\partial a}& (a-1)\frac{\partial w}{\partial b} \end{array}\right ).
\end{align*}
Next, we abelianize, setting $a$ equal to $x$ and $b$ equal to $y$, to obtain the Alexander matrix
$$
M= \left (
\begin{array}{cc} 
1-x^{\epsilon_2+\epsilon_4+\dots +\epsilon_{q-2}}y^{\epsilon_1+\epsilon_3+\dots +\epsilon_{q-1}}+(x-1)(\frac{\partial w}{\partial a})^\theta&(y-1)(\frac{\partial w}{\partial b})^\theta 
\end{array}
\right ).
$$
It remains to compute the partials $\frac{\partial w}{\partial a}$ and $\frac{\partial w}{\partial b}$ and their abelianizations  $(\frac{\partial w}{\partial a})^\theta$ and  $(\frac{\partial w}{\partial b})^\theta$.   

First note that if $\epsilon=\pm1$, then $\frac{\partial a^\epsilon}{\partial a}=\epsilon a^{(\epsilon-1)/2}$.  Proceeding by induction on $q$, we obtain
$$\frac{\partial w}{\partial a}=\sum_{i=1}^\frac{q-2}{2} \epsilon_{2 i}b^{\epsilon_1} a^{\epsilon_2} \dots b^{\epsilon_{2i-1}} a^\frac{\epsilon_{2i}-1}{2}\text{\quad and \quad}
\frac{\partial w}{\partial b}=\sum_{i=1}^\frac{q}{2} \epsilon_{2 i-1}b^{\epsilon_1} a^{\epsilon_2} \dots a^{\epsilon_{2i-2}} b^\frac{\epsilon_{2i-1}-1}{2}$$
which abelianize to
$$\left (\frac{\partial w}{\partial a}\right)^\theta=\sum_{i=1}^\frac{q-2}{2} \epsilon_{2i}x^{\frac{\epsilon_{2i}-1}{2}+\sum_{j=1}^{i-1}\epsilon_{2j}}y^{\sum_{k=1}^i \epsilon_{2k-1}}\text{ and }
\left (\frac{\partial w}{\partial b}\right)^\theta=\sum_{i=1}^\frac{q}{2} \epsilon_{2i-1}x^{\sum_{j=1}^{i-1}\epsilon_{2j}}y^{\frac{\epsilon_{2i-1}-1}{2}+\sum_{k=1}^{i-1} \epsilon_{2k-1}}
,
$$ respectively.
Now, using the fact that $(x-1)\epsilon t^{(\epsilon-1)/2}=x^\epsilon-1$, if $\epsilon=\pm1$, we obtain
\begin{equation}\label{first}(x-1)\left (\frac{\partial w}{\partial a}\right)^\theta=\sum_{i=1}^\frac{q-2}{2} (x^{\epsilon_{2 i}}-1)x^{\sum_{j=1}^{i-1}\epsilon_{2j}} y^{\sum_{k=1}^i\epsilon_{2k-1}}
\end{equation}
and similarly, 
\begin{equation}\label{second}(y-1)\left (\frac{\partial w}{\partial b}\right)^\theta=\sum_{i=1}^\frac{q}{2} (y^{\epsilon_{2 i}}-1)x^{\sum_{j=1}^{i-1}\epsilon_{2j}} y^{\sum_{k=1}^{i-1}\epsilon_{2k-1}}.
\end{equation}
Adding \eqref{first} and \eqref{second} gives
$$(x-1)\left (\frac{\partial w}{\partial a}\right)^\theta+(y-1)\left (\frac{\partial w}{\partial b}\right)^\theta=
-1+x^{\epsilon_2+\epsilon_4+\dots+\epsilon_{q-2}}y^{\epsilon_1+\epsilon_3+\dots+\epsilon_{q-1}}.$$
Hence the $(1,1)$-entry of the Alexander matrix is
$$M_{1,1}=-1+x^{\epsilon_2+\epsilon_4+\dots+\epsilon_{q-2}}y^{\epsilon_1+\epsilon_3+\dots+\epsilon_{q-1}}+(x-1)\left (\frac{\partial w}{\partial a}\right)^\theta=-(y-1)\left (\frac{\partial w}{\partial b}\right)^\theta.$$
Thus, the Alexander matrix is
$$M=\left(\begin{array}{cc} -(y-1)\left (\frac{\partial w}{\partial b}\right)^\theta&(x-1)\left (\frac{\partial w}{\partial b}\right)^\theta 
\end{array} 
\right ).
$$
The Alexander polynomial is the greatest common divisor of the two entries of $M$. Thus
$\Delta(x,y)=\left (\frac{\partial w}{\partial b}\right)^\theta$, giving the desired formula for $\Delta(x,y)$.
\hfill $\square$

\begin{corollary} Algorithm 3 follows from Equation~\eqref{2variableFormula}.
\end{corollary}
\noindent{\bf Proof:} For any $p$ and $q$, we have that $\epsilon_1=1$. Thus the first term of \eqref{2variableFormula} is 1. This corresponds
to starting the walk at the origin. Assume that the $i$-th term of \eqref{2variableFormula} is $\epsilon_{2i-1}x^ny^m$. The next
term is now $\epsilon_{2i-1}x^{n+\epsilon_{2i}}y^{m+(\epsilon_{2i-1}+\epsilon_{2i+1})/2}$. Thus we step right or left if $\epsilon_{2i}$ is positive or negative, 
respectively, and go up, stay level, or go down if $\epsilon_{2i-1}+\epsilon_{2i+1}$ is $2, 0$, or $-2$, respectively. Finally, it is easy to show by induction, that $\epsilon_{2i-1}=(-1)^{n+m}$. Thus the terms alternate sign in checkerboard fashion as claimed. \hfill $\square$



\end{document}